\documentstyle[amssymb]{amsart}

\numberwithin{equation}{section}

\theoremstyle{plain}
\newtheorem{thm}{Theorem}[section]
\newtheorem{prop}[thm]{Proposition}
\newtheorem{lem}[thm]{Lemma}

\theoremstyle{definition}

\theoremstyle{remark}
\newtheorem{rem}[thm]{Remark}

\newcommand{\comment}[1]{}
\newcommand{\Label}[1]{\label{#1}}
\newcommand{\br}[1]{\langle #1 \rangle}
\newcommand{\Pbr}[1]{\big\{{ #1 }\big\}}

\newcommand{\BZ}{{\Bbb Z}}

\newcommand{\BC}{{\Bbb C}}

\title[\ ]{
HIgher order Painlev\'e equations of type $A^{(1)}_l$ 
}

\begin{document}
%\begin{flushright}
%{\small {\em Preliminary Version}: \today}
%\end{flushright}

\par\bigskip
\maketitle
\begin{center}
Masatoshi Noumi  
and 
Yasuhiko Yamada 
\end{center}

{\small
\begin{center}
Department of Mathematics, Kobe University, 
Rokko, Kobe 657-8501, Japan
\end{center}
}

\begin{abstract}
A series of systems of nonlinear equations with affine Weyl group symmetry 
of type $A^{(1)}_l$ is studied.  
This series gives a generalization of Painlev\'e equations $P_{\,\text{IV}}$
and $P_{\,\text{V}}$ to higher orders. 
\end{abstract}

%%%%%%%%%%
\section*{Introduction}
In this paper we propose a series of systems of nonlinear differential equations 
which have symmetry under the affine Weyl groups of type $A^{(1)}_l$ ($l\ge 2$). 
These systems are thought of as higher order analogues of the Painlev\'e equations 
$P_{\,\text{IV}}$ and $P_{\,\text{V}}$.  
We will study in particular their Hamiltonian structures and 
some basic properties of their $\tau$-functions. 

\medskip
For each $l=2,3,\ldots$, we consider a differential system for 
$(l+1)$ unknown functions $f_0,\ldots,f_l$, containing 
complex parameters $\alpha_0,\ldots,\alpha_l$ 
which correspond to the 
simple roots of the affine root system of type $A^{(1)}_l$.  
In what follows, we set $\alpha_0+\cdots+\alpha_l=k$. 
Our differential system is defined as follows according to the parity of $l$: 
For $l=2n$ ($n=1,2,\ldots$)
\begin{equation}\Label{A2n}
A^{(1)}_{2n}:\qquad f_j'=f_j\big (\sum_{1\le r\le n}f_{j+2r-1}
-\sum_{1\le r\le n} f_{j+2r}\big)+\alpha_j
\end{equation}
where $0\le j \le 2n$, 
and for $l=2n+1$ ($n=1,2,\ldots$)
\begin{equation}\Label{A2n+1}
\begin{align}
A^{(1)}_{2n+1}:\quad f_j'&=f_j\big(
\sum_{ 1\le r\le s\le n } f_{j+2r-1}f_{j+2s}
-\sum_{ 1\le r\le s\le n} f_{j+2r}f_{j+2s+1}
\big)
\\ 
&\quad+\big(\frac{k}{2}-\sum_{1\le r\le n}\alpha_{j+2r}\big)f_j +
\alpha_j\big(\sum_{1\le r\le n}f_{j+2r}\big),\notag
\end{align}
\end{equation}
where $0\le j\le 2n+1$. 
In \eqref{A2n} and \eqref{A2n+1}, 
$'$ stands for the derivation $d/dt$ with respect to 
an independent variable $t$.  
A table of formulas of $f_0'$ for small $l$ is given in Appendix for convenience. 
Formulas for the other $f_j'$ are obtained from those simply by the rotation of indices.
It can be shown that the differential systems for $l=2$ and $l=3$ are equivalent to 
the fourth and the fifth Painlev\'e equations, respectively(see \cite{NY5}).

This paper is organized as follows.  In Section 1, we discuss symmetry of our system 
under the affine Weyl group of type $A^{(1)}_l$, by describing explicit B\"acklund 
transformations.  
After formulating a Poisson structure for our system, we will 
construct in Section 3 
certain canonical coordinates with respect to the Poisson structure, 
and show that our system can be equivalently written as a Hamiltonian 
system with a polynomial Hamiltonian.  
In the final section, we introduce a family of $\tau$-functions for our system 
and discuss their B\"acklund transformations.  In particular 
we will show that the variables $f_j$ are expressed multiplicatively  
in terms of the $\tau$-functions and their B\"acklund transformations. 
This shows that our $\tau$-functions are consistent with those introduced 
in \cite{NY3} from the viewpoint of discrete dynamical systems. 
%%%%%%%%%%%%%%%%%%%
\section{Affine Weyl group symmetry} 
Let $\BC(\alpha;f)$ be the field of rational functions in 
$\alpha=(\alpha_0,\ldots,\alpha_l)$ and $f=(f_0,\ldots,f_l$).  
Then the differential system \eqref{A2n} (resp. \eqref{A2n+1}) 
defines the structure of a differential field on $\BC(\alpha;f)$. 
We say that an automorphism of $\BC(\alpha;f)$ is a 
{\em B\"acklund transformation} of the differential system 
if it commutes with the derivation $'$.  

For each $i=0,1,\ldots,l$, we define an automorphism $s_i$ of $\BC(\alpha;f)$
as follows:
\begin{equation}\Label{WAl1}
\begin{array}{lllll}
s_i(\alpha_i)=-\alpha_i, 
&s_i(\alpha_{j})=\alpha_j+\alpha_i&(j=i\pm 1),
&s_i(\alpha_{j})=\alpha_j&(j\ne i,i\pm1 ), \\
s_i(f_i)=f_i,
&s_i(f_j)=\displaystyle{f_j\pm\frac{\alpha_i}{f_i}}&(j=i\pm 1),
&s_i(f_j)=f_j&(j\ne i,i\pm1),
\end{array}
\end{equation}
where the indices $0,1,\ldots,l$ are understood as elements of $\BZ/(l+1)\BZ$. 
We also define an automorphism $\pi$ of $\BC(\alpha;f)$ by 
\begin{equation}\Label{WAl2}
\pi(\alpha_j)=\alpha_{j+1},\quad \pi(f_j)=f_{j+1}.
\end{equation}

\begin{thm}[\cite{NY3}]
The automorphisms $s_0,s_1,\ldots,s_l$ and $\pi$ described above define 
a representation of the extended affine Weyl group 
$\widetilde{W}=\br{s_0,\ldots,s_l,\pi}$ of type $A^{(1)}_{l}$. 
Namely, they satisfy the commutation relations
\begin{equation}\Label{rel-s}
s_i^2=1, \quad s_is_j=s_js_i\ \ (j\ne i,i\pm1),\quad
s_is_js_i=s_js_is_j \ \ (j=i\pm1)
\end{equation}
for $i,j=0,1,\ldots,l$ and 
\begin{equation}
\pi^{l+1}=1,\qquad \pi s_{i}=s_{i+1} \pi\quad(i=0,1,\ldots,l).
\end{equation}
\end{thm}
\noindent
Note that the extended affine Weyl group 
$\widetilde{W}=W\rtimes\{1,\pi,\ldots,\pi^l\}$ is the extension of 
the ordinary affine Weyl group $W=\br{s_0,\ldots,s_l}$ of type $A^{(1)}_l$ 
by a cyclic group of order $l+1$ generated by the {\em diagram rotation} $\pi$.

\begin{thm}\Label{thmA}
The action of $\widetilde{W}$ defined as 
above commutes with the derivation of the differential field 
$\BC(\alpha;f)$.
\end{thm}

\noindent
In this sense, our differential system \eqref{A2n} (resp. \eqref{A2n+1}) 
admits the action of the extended affine Weyl group $\widetilde{W}$ of 
type $A^{(1)}_l$, as a group of B\"acklund transformations. 
Theorem \ref{thmA} can be checked essentially by direct computations.  We will explain below 
how such computations can be carried out. 

\begin{rem}
Let $A=(a_{ij})_{0\le i,j\le l}$ be the generalized Cartan matrix of type 
$A^{(1)}_{l}$:
\begin{equation}\Label{GCM}
a_{jj}=2,\quad a_{ij}=-1\ \ (j=i\pm1),\quad a_{ij}=0\ \ (j\ne i,i\pm1). 
\end{equation}
Then the action of $s_0,\ldots,s_l$ on the simple roots $\alpha_0,\ldots,\alpha_l$ 
is described as 
\begin{equation}\Label{s-on-a}
s_i(\alpha_j)=\alpha_j-\alpha_i a_{ij}\quad(i,j=0,1,\ldots,l). 
\end{equation}
We introduce an $(l+1)\times(l+1)$ matrix $U=(u_{ij})_{0\le i,j\le l}$ by setting
\begin{equation}\Label{def-U}
u_{ij}=\pm 1 \quad (j=i\pm 1), \quad u_{ij}=0\quad(j\ne i\pm 1). 
\end{equation}
Then the B\"acklund transformations $s_0,\ldots,s_l$ are determined by the 
formulas
\begin{equation}\Label{s-on-f}
s_{i}(f_j)=f_j+\frac{\alpha_i}{f_i} u_{ij}\quad(i,j=0,1,\ldots,l). 
\end{equation}
For a general treatment of affine Weyl group symmetry in terms of these data 
$A$ and $U$, as well as discrete dynamical systems arising from B\"acklund 
transformations, we refer the reader to \cite{NY3}. 
\end{rem}

\medskip
For practical computations of B\"acklund transformations, 
it is convenient to use the {\em Demazure operators\/}
$\Delta_i$ $(i=0,\ldots,l)$ defined by
\begin{equation}
\Delta_i(\varphi)=\frac{1}{\alpha_i}
(s_i(\varphi)-\varphi)\qquad (\varphi\in \BC(\alpha;f)). 
\end{equation}
The action of $\Delta_i$ on a product can be determined 
by the {\em twisted\/} Leibniz rule
\begin{equation}
\Delta_i(\varphi\psi)=\Delta_i(\varphi) \psi+s_i(\varphi) \Delta_i(\psi).  
\end{equation}
We remark that $\Delta_i(\varphi)=0$ if and only if $\varphi$ is 
$s_i$-invariant and that one has 
$\Delta_i(\varphi\psi)=\varphi \Delta(\psi)$ 
when $\varphi$ is $s_i$-invariant.
Note that \eqref{s-on-a} and \eqref{s-on-f} can be rewritten as 
\begin{equation}
\Delta_i(\alpha_j)=-a_{ij},\quad \Delta_i(f_j)=\frac{u_{ij}}{f_i}
\quad(i,j=0,\ldots,l), 
\end{equation}
in terms of the Demazure operators. 
Namely one has
\begin{equation}\Label{D-on-sf}
\begin{array}{lllll}
\Delta_i(\alpha_i)=-2, 
&\Delta_i(\alpha_{j})=1&(j=i\pm 1),
&\Delta_i(\alpha_{j})=0&(j\ne i,i\pm1 ),\\
\Delta_i(f_i)=0,
&\Delta_i(f_j)=\displaystyle{\frac{\pm1}{f_i}}&(j=i\pm 1),
&\Delta_i(f_j)=0&(j\ne i,i\pm1). 
\end{array}
\end{equation}
It is also well known that the relations 
\eqref{rel-s} imply 
\begin{equation}
\Delta_i^2=0,\ \Delta_i \Delta_j =\Delta_j \Delta_i \ (j\ne i\pm 1),\ \ 
\Delta_i \Delta_j \Delta_i=\Delta_j \Delta_i \Delta_j \  (j=i\pm 1), 
\end{equation}
for $i,j=0,\ldots,l$. 

\medskip\noindent
{\em Proof of Theorem \ref{thmA}}.
For each $j=0,\ldots,l$, we denote by $F_j=F_j(\alpha;f)$ the polynomial 
appearing on the right-hand side of \eqref{A2n} or \eqref{A2n+1}.
Then it is easy to see that $s_i(f_j)'=s_i(f_j')$ if and only if 
\begin{equation}
F_j-\frac{\alpha_i}{f_i^2}\,F_i \,  u_{ij}=s_i(F_j)\qquad(i,j=0,\ldots,l)
\end{equation}
or equivalently, 
\begin{equation}\Label{D-on-F}
\Delta_i(F_j)=-\frac{u_{ij}}{f_i^2} \, F_i\qquad(i,j=0,\ldots,l)
\end{equation}
in terms of the Demazure operators. 
By the rotation symmetry, it is enough to check \eqref{D-on-F} for $j=0$:
\begin{equation}\Label{tochk}
\Delta_{1}(F_0)=\frac{F_1}{f_1^2}, \quad\Delta_{l}(F_0)=-\frac{F_l}{f_l^2},
\quad \Delta_i(F_0)=0\ \ (i\ne 0,l). 
\end{equation}
We will show for example the equality $\Delta_{1}(F_0)=F_1/f_1^2$ for 
$l=2n$ and $l=2n+1$, separately. 
When $l=2n$, we have
\begin{equation}
F_0=f_0(f_1-f_2+\cdots-f_{2n})+\alpha_0. 
\end{equation}
Since $\Delta_1(f_j)=0$ for $j\ne 0,2$, we have
\begin{eqnarray}
\Delta_1(F_0)&=&\Delta_1(f_0)(\sum_{i=1}^{2n}(-1)^{i-1} f_i) 
-s_1(f_0)\Delta_1(f_2)+\Delta_1(\alpha_0)\\
&=&-\frac{1}{f_1}(\sum_{i=1}^{2n}(-1)^{i-1} f_i)
-(f_0-\frac{\alpha_1}{f_1})\frac{1}{f_1}+1
\notag\\
&=&\frac{1}{f_1}(\sum_{i=2}^{2n}(-1)^{i} f_i -f_0)+\frac{\alpha_1}{f_1^2}
=\frac{1}{f_1^2}\, F_1. 
\notag
\end{eqnarray}
When $l=2n+1$, we have 
\begin{equation}
\begin{align}
F_0=&f_0(f_1f_2+f_1f_4+\cdots+f_{2n-1}f_{2n}-
f_2f_3-\cdots-f_{2n}f_{2n+1})\\
&
+(\frac{k}{2}-\alpha_2-\alpha_4-\cdots-\alpha_{2n})f_0+\alpha_0(f_2+f_4+\cdots+f_{2n}).
\notag
\end{align}
\end{equation}
Hence, $\Delta_1(F_0)$ is computed as follows:
\allowdisplaybreaks
\begin{eqnarray*}
%&\Delta_1(F_0)\notag\\
&&\Delta_1(f_0)
\big(
\sum_{1\le r\le s \le n} f_{2r-1}f_{2s}-\sum_{1\le r\le s\le n} f_{2r}f_{2s+1}\big)
%%\\
%%&\quad
+s_1(f_0)\Delta_1(f_2)\big(f_1-\sum_{r=1}^n f_{2r+1}\big)
%\notag
\\
&&
-\Delta_1(\alpha_{2})f_0+s_1(\frac{k}{2}-\sum_{r=1}^n \alpha_{2r})\Delta_1(f_0)
%%\notag\\
%%&\quad
+\Delta_1(\alpha_0)\big(\sum_{r=1}^n f_{2r}\big)
+s_1(\alpha_0)\Delta_1(f_2)
\notag\\
&=&-\frac{1}{f_1}
\big(
\sum_{1\le r\le s \le n} f_{2r-1}f_{2s}-\sum_{1\le r\le s\le n} f_{2r}f_{2s+1}\big)
%\notag\\
%&\quad
+(f_0-\frac{\alpha_1}{f_1})\frac{1}{f_1}\big(f_1-\sum_{r=1}^n f_{2r+1}\big)
\notag\\
&&\quad
-f_0+(\frac{k}{2}-\sum_{r=1}^n \alpha_{2r+1}-\alpha_0)\frac{1}{f_1}
%%\notag\\
%%&\quad
+\big(\sum_{r=1}^n f_{2r}\big)
+(\alpha_0+\alpha_1)\frac{1}{f_1}
\notag\\
&=&\frac{1}{f_1}
\big(\sum_{1\le r\le s\le n} f_{2r}f_{2s+1}
-\sum_{2\le r\le s \le n} f_{2r-1}f_{2s}-\sum_{r=1}^n f_{2r+1}f_0\big)
\notag\\
&&\quad
+(\frac{k}{2}-\sum_{r=1}^n \alpha_{2r+1})\frac{1}{f_1}
+\frac{\alpha_1}{f_1^2}\big(\sum_{r=1}^n f_{2r+1}\big)
=\frac{1}{f_1^2}F_1.
\notag
\end{eqnarray*}
The other formulas in \eqref{tochk} can be checked similarly.
\quad\qed

\section{Poisson structure}
By using the matrix $U$ defined in \eqref{def-U}, we introduce the {\em Poisson bracket} 
$\Pbr{\ ,\ }$ on $\BC(\alpha;f)$ as follows: 
\begin{equation}\Label{def-Pbr}
\Pbr{\varphi,\psi} = \sum_{0\le i,j \le l}
\frac{\partial \varphi}{\partial f_i}\ u_{ij}\ \frac{\partial \psi}{\partial f_j}. 
\end{equation}
Note that $\Pbr{\ ,\ }$ is a $\BC(\alpha)$-bilinear skewsymmetric form such that
\begin{equation}
\Pbr{f_i,f_j}=u_{ij}\quad(i,j=0,1,\ldots,l),
\end{equation}
namely,
\begin{equation}\Label{Pbr-on-f}
\Pbr{f_i,f_j}=\pm 1\quad(j=i\pm1),\quad \Pbr{f_i,f_j}=0\quad(j\ne i\pm 1)
\end{equation}
for $i,j=0,1,\ldots,l$.
By direct calculations, one can show
\begin{prop} 
The skewsymmetric form $\Pbr{\ ,\ }$ above defines a $\widetilde{W}$-invariant 
Poisson structure on the differential field $\BC(\alpha;f)$. 
Namely, one has
\begin{equation*}
\begin{align*}
(1)&\ \ 
\Pbr{\varphi,\psi_1\psi_2}=
\psi_1\Pbr{\varphi,\psi_2}+\Pbr{\varphi,\psi_1}\psi_2,
\Pbr{\varphi_1\varphi_2,\psi}=
\varphi_1\Pbr{\varphi_2,\psi}+\Pbr{\varphi_1,\psi}\varphi_2,\\
(2)&\ \ 
\Pbr{\varphi_1,\Pbr{\varphi_2,\varphi_3}}+
\Pbr{\varphi_2,\Pbr{\varphi_3,\varphi_1}}+
\Pbr{\varphi_3,\Pbr{\varphi_1,\varphi_2}}=0,\\
(3)&\ \ 
w(\Pbr{\varphi,\psi})=\Pbr{w(\varphi),w(\psi)} \quad
\text{for any}\ \  \varphi,\psi\in \BC(\alpha;f) \ \ \text{and} \ \ w\in \widetilde{W}. 
\end{align*}
\end{equation*}
\end{prop}

We remark that our Poisson bracket has a nontrivial radical. 
Regard $\Pbr{\ ,\ }$ as a skewsymmetric form on the $(l+1)$-dimensional 
vector space $E=\bigoplus_{j=0}^l \BC(\alpha)f_j$. 
Then its radical is precisely the subspace of all linear combinations 
$c_0f_0+\cdots+c_l f_l$ with coefficients in $\BC(\alpha)$  
such that $\sum_{j=0}^l u_{ij}c_j=0 $ for all $i=0,\ldots,l$.  
It is a one-dimensional subspace generated by
\begin{equation}\Label{def-g}
g=f_0+f_1+\cdots+f_{2n}
\end{equation}
if $l=2n$, and is a two-dimensional subspace generated by
\begin{equation}\Label{def-g01}
g_0=f_0+f_2+\cdots+f_{2n},\quad g_1=f_1+f_3+\cdots+f_{2n+1},
\end{equation}
if $l=2n+1$, respectively. 
Note that $g$ is $\widetilde{W}$-invariant when $l=2n$, and that
$g_0,g_1$ are $W$-invariant, and that $g_0+g_1$ is $\widetilde{W}$-invariant
when $l=2n+1$.

\medskip
We will describe below our differential systems \eqref{A2n} and \eqref{A2n+1}
by means of the Poisson structure introduced above. 
In order to define a ``Hamiltonian'' we need to fix some notation. 

\medskip
For each $i=1,\ldots,l$, we denote by $\varpi_i$  
the $i$-th {\em fundamental weight\/} of the 
finite root system of type $A_l$ 
\begin{eqnarray}\Label{def-fwt}
\varpi_i&=&\frac{1}{l+1}\big\{
(l+1-i)\sum_{r=1}^i r\alpha_r + i\sum_{r=i+1}^{l}(l+1-r)\alpha_r\big\}\\
&=&\sum_{r=1}^l\big(\min\{i,r\}-\frac{ir}{l+1}\big)\alpha_r \notag
\end{eqnarray} 
and set $\varpi_0=0$.  
Note that $\varpi_1,\ldots,\varpi_l$ form the dual basis of the simple roots 
$\alpha_1,\ldots,\alpha_l$ with respect to the symmetric bilinear form 
$\br{\ ,\ }$ such that $\br{\alpha_i,\alpha_j}=a_{ij}$ ($i,j=1,\ldots,l$). 
The simple affine roots $\alpha_0,\ldots,\alpha_l$ 
are then expressed  as 
\begin{equation}
\alpha_j=-\varpi_{j-1}+2 \varpi_j -\varpi_{j+1}+\delta_{j,0} k \quad(j=1,\ldots,l), 
\end{equation} 
in terms of the fundamental weights $\varpi_1,\ldots,\varpi_l$.  
We remark that the reflections $s_0,\ldots,s_l$ act on the fundamental weights 
as 
\begin{equation}
s_0(\varpi_j)=\varpi_j+\alpha_0,\quad
s_i(\varpi_j)=\varpi_j-\delta_{i,j} \alpha_i \ \ (i=1,\ldots,l),
\quad 
\end{equation}
for each $j=1,\ldots,l$. 
Notice that the diagram rotation $\pi$ acts on $\varpi_1,\ldots,\varpi_l$ 
nontrivially. 

Let $\Gamma$ be the Dynkin diagram of type $A^{(1)}_l$; it is a circle with 
$(l+1)$ nodes labeled by the elements of $\BZ/(l+1)\BZ$. 
For each chain $C$ of $\Gamma$, consisting of
consecutive nodes $j,j+1,\ldots,j+m-1$ ($m\le l$), we denote by 
$\chi(C)$ the alternating sum of corresponding fundamental weights: 
\begin{equation}\Label{chi}
\chi(C)=\varpi_j-\varpi_{j+1}+\ldots+(-1)^{m-1}\varpi_{j+m-1}. 
\end{equation} 
For each subdiagram $C$ of $\Gamma$ with $C\ne \Gamma$, 
%all the conneted components are chains,  
we denote by $\chi(C_i)=\sum_i \chi(C_i)$ the sum of $\chi(C_i)$ over all 
connected components $C_i$. 

For each $d=1,\ldots,l+1$, we denote by ${\cal S}_d$ the set of all subsets 
$K\subset\{0,1,\ldots,l\}$ with cardinality $d$ such that the connected components 
of the diagram $K^{\text{c}}=\Gamma\backslash K$,
obtained by removing the nodes of $K$, are all chains of
{\em even} nodes. 
For each $K\in {\cal S}_d$, we set $f_K=\prod_{i\in K} f_i$. 
We remark that the set ${\cal S}_d$ is nonempty if and only if $l+1-d$ is even, 
and that a subset $K$ of $\{0,\ldots,l\}$ with $|K|=d$ then belongs to ${\cal S}_d$ 
if and only if it has an expression $K=\{k_1,k_2,\ldots,k_d\}$ with a sequence 
$0\le k_1<k_2<\cdots<k_d\le l$ such that
\begin{equation}
(k_1,k_2,\ldots,k_d)\equiv (0,1,0,\ldots) \ \ \text{or} \ \ (1,0,1,\ldots)\mod 2. 
\end{equation}
With this expression, $\chi(K^{\text{c}})$ can be written as follows:
\begin{equation}
\chi(K^{\text{c}})=\sum_{i=0}^{d-1}\ \sum_{r=1}^{k_{i+1}-k_i-1}(-1)^{r-1}\varpi_{k_i+r}
\end{equation}
where $k_0=k_d-l-1$. 

\medskip
With the notation as above, 
we define a {\em Hamiltonian} $h_0$ for our system by
\begin{equation}\Label{h0-2n}
\begin{align}
h_0&=\sum_{K\in{\cal S}_3} f_K +\sum_{K\in{\cal S}_1} \chi(K^{\text{c}})\, f_K\\
&=\sum f_{k_1} f_{k_2} f_{k_3}   
+\sum_{i=0}^{2n} \big(\sum_{r=1}^{2n}(-1)^{r-1}\varpi_{i+r}\big)\,f_i
\notag
\end{align}
\end{equation}
when $l=2n$, 
where the summation $\sum f_{k_1} f_{k_2} f_{k_3}$ is taken over all triples 
$(k_1,k_2,k_3)$ such that $0\le k_1<k_2<k_3\le 2n$ and that 
$(k_1,k_2,k_3)\equiv (0,1,0)$ or $(1,0,1) \mod 2$.  
%We remark that $|{\cal S}_3|=n(n+1)(2n+1)/6$ in this case.
When $l=2n+1$, we define $h_0$ by
\begin{equation}\Label{h0-2n+1}
h_0=\sum_{K\in{\cal S}_4} f_K + \sum_{K\in{\cal S}_2} \chi(K^{\text{c}})\,f_K
+\big(\sum_{i=1}^{2n+1}(-1)^{i-1} \varpi_i\big)^2.
\end{equation}
We remark that, when $l=2n+1$, the sum appearing in the constant term 
of \eqref{h0-2n+1} has an alternative 
expression 
\begin{equation}
\sum_{i=1}^{2n+1}(-1)^{i-1} \varpi_i=\frac{1}{2}\sum_{r=0}^n \alpha_{2r+1}
\end{equation}
by simple roots. 
We give in Appendix some explicit formulas of $h_0$ for small $l$ for convenience. 

\begin{prop} \Label{prop-by-Pbr}
When $l=2n$, the differential system \eqref{A2n} can be expressed as follows:
\begin{equation}\Label{A2n-f-by-Pbr}
f_j'=\Pbr{h_0,f_j}+\delta_{j,0}\, k\quad(j=0,\ldots,2n). 
\end{equation}
Hence one has 
\begin{equation}\Label{A2n-by-Pbr}
\varphi'=\Pbr{h_0,\varphi}+k \frac{\partial \varphi}{\partial f_0}
\quad\text{for any}\quad \varphi\in \BC(\alpha;f). 
\end{equation}
When $l=2n+1$, the differential system \eqref{A2n+1} can be expressed as follows:
\begin{equation}\Label{A2n+1-f-by-Pbr}
f_j'=\Pbr{h_0,f_j}-(-1)^{j}\, \frac{k}{2}\,f_j+\delta_{j,0}\, k\, g_0
\quad(j=0,\ldots,2n+1),
\end{equation}
with $g_0$ defined by \eqref{def-g01}.
Hence one has 
\begin{equation}\Label{A2n+1-by-Pbr}
\varphi'=\Pbr{h_0,\varphi}-
\frac{k}{2}
\big(\sum_{i=0}^{2n+1}(-1)^i f_i\frac{\partial\varphi }{\partial f_i}\big)
+ k g_0 \frac{\partial \varphi}{\partial f_0}
\end{equation} 
for any $\varphi\in \BC(\alpha;f)$.
\end{prop}
\begin{pf}
%Proposition \ref{prop-by-Pbr} can be proved by a combinatorial argument 
%from the definition of $h_0$, \eqref{h0-2n} and \eqref{h0-2n+1}. 
For each $j=0,\ldots,l$, we define the vector field $X_j$ by
\begin{equation}
X_j(\varphi)=\Pbr{\varphi,f_j}
=\big(\frac{\partial\quad}{\partial f_{j-1}}-\frac{\partial\quad}{\partial f_{j+1}}\big)
\varphi. 
\end{equation}
We now consider the case of $A^{(1)}_{2n}$.
From the definition \eqref{h0-2n}, $X_j(h_0)$ is computed as follows:
\allowdisplaybreaks
\begin{eqnarray*}
&&\sum_{K\in {\cal S}_2(\Gamma\backslash\{j-1\})} f_K
-\sum_{K\in {\cal S}_2(\Gamma\backslash\{j+1\})} f_K
%\\
%&\qquad
+\chi(\Gamma\backslash\{j-1\})-\chi(\Gamma\backslash\{j+1\})
\notag\\
&=&f_j\sum_{r=1}^{n}f_{j+2r-1}-f_{j}\sum_{r=1}^n f_{j+2r}
-\varpi_{j-1}+2\varpi_j-\varpi_{j+1}
\notag
\\
&=&f_j\big(\sum_{i=1}^{2n} (-1)^{i-1}f_{j+i}\big)+\alpha_j-\delta_{j,0}\, k
=F_j-\delta_{j,0}\, k,
\notag
\end{eqnarray*} 
with the notation of ${\cal S}_d$ extended to subdiagrams of 
$\Gamma=\{0,1,\ldots,l\}$.  This proves \eqref{A2n-f-by-Pbr}.
(For a subdiagram $G$ of $\Gamma$, ${\cal S}_d(G)$ stands for the set of all 
subsets $K\subset G$ with $|K|=d$ such that the connected components of $G\backslash K$
are all chains of even nodes.)
Formulas \eqref{A2n+1-f-by-Pbr} for the case of $A^{(1)}_{2n+1}$ can 
be established in a similar way. 
In fact $X_j(h_0)$ is given by 
\allowdisplaybreaks
\begin{eqnarray*}
\allowdisplaybreaks
%X_j(h_0)
&&\sum_{K\in{\cal S}_3(\Gamma\backslash\{j-1\})} f_K
-\sum_{K\in{\cal S}_3(\Gamma\backslash\{j+1\})} f_K
\notag\\*
&&\quad
+\sum_{\{i\}\in{\cal S}_1(\Gamma\backslash\{j-1\})} \chi(\Gamma\backslash\{j-1,i\})f_i
-\sum_{\{i\}\in{\cal S}_1(\Gamma\backslash\{j+1\})} \chi(\Gamma\backslash\{j+1,i\})f_i
\notag\\
&=&f_j\big(\sum_{K\in{\cal S}_2([j+1,j-2])}f_K
-\sum_{K\in{\cal S}_2([j+2,j-1])}f_K
\big)
\notag\\
&&\quad
+\big(-\varpi_{j-1}+2\varpi_j-\varpi_{j+1}+
 (-1)^j 2\sum_{i=1}^{2n+1}(-1)^{i-1}\varpi_{i}\big)\, f_j
\notag\\
&&\quad
+(-\varpi_{j-1}+2\varpi_j-\varpi_{j+1})\sum_{r=1}^{n}\,f_{j+2r}
\notag\\
&=&f_j\big(\sum_{K\in{\cal S}_2([j+1,j-2])}f_K
-\sum_{K\in{\cal S}_2([j+2,j-1])}f_K
\big)
\notag\\
&&\quad
+\big(\alpha_j-\delta_{j,0}\, k +
 (-1)^j \sum_{r=0}^{n}\alpha_{2r+1}\big)\, f_j
+(\alpha_j-\delta_{j,0}\, k)\sum_{r=1}^{n}\,f_{j+2r}
\notag\\
&=&F_j+(-1)^j \frac{k}{2}f_j-\delta_{j,0}\,k\,\big(\sum_{r=0}^n f_{2r}\big),
\notag
\end{eqnarray*}
where $[a,b]$ stands for the chain $\{a,a+1,\ldots,b-1,b\}$. 
This proves \eqref{A2n+1-f-by-Pbr}.
\end{pf}

\section{Canonical coordinates and the Hamiltonian system}
By using Proposition \ref{prop-by-Pbr}, we introduce {\em canonical coordinates}
for our differential system. 
We discuss the two cases of $A^{(1)}_{2n}$ and $A^{(1)}_{2n+1}$ separately.

\par\medskip\noindent
{\bf Case $A^{(1)}_{2n}$: }\quad
In view of \eqref{Pbr-on-f},
we define a new coordinate system 
\begin{equation}
(q;p;x)=(q_1,\dots,q_n;p_1,\ldots,p_n;x)
\end{equation}
for the affine space with coordinates $f=(f_0,\ldots,f_{2n})$ as 
follows:
\begin{equation}
\begin{align}
&q_1=f_2, \ \ q_2=f_4, \ \  \ldots, \ \ q_n=f_{2n},\\ 
&p_1=f_1, \ \ p_2=f_1+f_3,\ \  \ldots,\ \  p_n=f_1+f_3+\cdots+f_{2n-1},\notag\\
&x=g=f_0+f_1+\cdots+f_{2n}. \notag
\end{align}
\end{equation} 
Note that the inverse of this coordinate transformation is given by
\begin{equation}\Label{f-in-qp-2n}
\begin{align}
&f_0=x-q_1-q_2-\cdots-q_n-p_n,\\
&f_1=p_1,\ \  f_2=q_1, \ \ f_3=p_2-p_1, \ \ f_4=q_2,\ldots,\notag\\
&f_{2n-1}=p_n-p_{n-1}, \ \ f_{2n}=q_n.\notag
\end{align}
\end{equation}
Then it is easy to show that 
\begin{equation}
\Pbr{p_i,q_j}=\delta_{i,j},\quad 
\Pbr{q_i,q_j}=\Pbr{p_i,p_j}=\Pbr{p_i,x}=\Pbr{q_i,x}=0 
\end{equation}
for $i,j=1,\ldots,n$.  
Hence we have
\begin{equation}\Label{Pbr-can}
\Pbr{\varphi,\psi}=\sum_{i,j=1}^n \bigg(
\frac{\partial\varphi}{\partial p_i}\frac{\partial \psi}{\partial q_i}
-\frac{\partial\varphi}{\partial q_i}\frac{\partial \psi}{\partial p_i}\bigg)
\quad (\varphi,\psi\in\BC(\alpha;f)). 
\end{equation}
Noting that $\partial/\partial f_0=\partial/\partial x$, 
from \eqref{A2n-by-Pbr} we obtain
\begin{equation}
\varphi'=\Pbr{h_0,\varphi} + k \frac{\partial\varphi}{\partial x}
\quad\text{for any}\quad \varphi\in\BC(\alpha;f).
\end{equation}
Let us denote by $H=H(q;p;x)\in\BC(\alpha)[q;p;x]$ the polynomial in $(q;p)$ which 
represents our Hamiltonian $h_0$ in the coordinates $(q;p;x)$. 
Then we see that the differential system \eqref{A2n} is equivalent to the 
{\em Hamiltonian system}
\begin{equation}
\frac{dq_i}{dt}=\Pbr{H,q_i}=\frac{\partial H}{\partial p_i},\quad
\frac{dp_i}{dt}=\Pbr{H,p_i}=-\frac{\partial H}{\partial q_i},\quad
\frac{dx}{dt}=k, 
\end{equation}
where $i=1,\ldots,n$. 
The Hamiltonian $H$ is determined explicitly as follow:
\begin{equation}\Label{H2n}
\begin{align}
H=& \big(x-\sum_{i=1}^n q_i \big) \big(\sum_{i=1}^n q_i p_i\big)
-\sum_{i=1}^n q_i p_i^2 
-\sum_{1\le i<j\le n}q_i(p_i-p_j)q_j \\
&
-\sum_{i=1}^n \big(\sum_{r=1}^i \alpha_{2r-1}\big) q_i+\sum_{i=1}^n \alpha_{2i} p_i
+\beta x,
\notag
%\\
%&
%+\frac{1}{2n+1}\big(\sum_{r=1}^n (n+1-r)\alpha_{2r-1}-r\alpha_{2r}\big)x.
%\notag
\end{align}
\end{equation}
where
\begin{equation}
\beta=\sum_{i=1}^{2n}(-1)^{i-1}\varpi_i=
\frac{1}{2n+1}\sum_{r=1}^n \big((n+1-r)\alpha_{2r-1}-r\alpha_{2r}\big).
\end{equation}

\medskip\noindent
{\bf Case $A^{(1)}_{2n+1}$:}\quad
Note first that \eqref{A2n+1-f-by-Pbr} implies
\begin{equation}
g_0'=\frac{k}{2} g_0,\quad g_1'=\frac{k}{2} g_1. 
\end{equation}
Hence, by setting
\begin{equation}
\widetilde{f}_{2r}=g_0 f_{2r}, \quad \widetilde{f}_{2r+1}=g_0^{-1} f_{2r+1}
\quad(r=0,1,\ldots,n), 
\end{equation}
we obtain
\begin{equation}
\widetilde{f}_j' = \Pbr{h_0,\widetilde{f}_j} + \delta_{j,0}\, k \, g_0^2
\quad(j=0,1,\ldots,2n+1).  
\end{equation}
We now introduce a new coordinate system 
\begin{equation}
(q;p;x)=(q_1,\ldots,q_n;p_1\ldots,p_n;x_0,x_1)
\end{equation}
 as follows:
\begin{equation}
\begin{align}
&q_1=g_0 f_2, \ \ q_2=g_0 f_4, \ \  \ldots, \ \ q_n=g_0 f_{2n},\\ 
&p_1=g_0^{-1}f_1, \ \ p_2=g_0^{-1}(f_1+f_3),\ \  \ldots,\ \  
p_n=g_0^{-1}(f_1+f_3+\cdots+f_{2n-1}),\notag\\
&x_0=g_0=f_0+f_2+\cdots+f_{2n},\ x_1=g_1=f_1+f_3+\cdots+f_{2n+1}. \notag
\end{align}
\end{equation}
The inverse transformation is then given by
\begin{equation}\Label{f-in-qp-2n+1}
\begin{align}
&f_0=x_0-x_0^{-1}(q_1+\cdots+q_n),\\
&f_1=x_0 p_1,\ \  f_2=x_0^{-1} q_1, \ \ f_3=x_0(p_2-p_1), 
\ \ f_4=x_0^{-1}q_2,\ldots,\notag\\
&f_{2n-1}=x_0(p_n-p_{n-1}), \ \ f_{2n}=x_0^{-1}q_n,
\ \ f_{2n+1}=x_1-x_0 p_n.\notag
\end{align}
\end{equation}
Since
\begin{equation}
\begin{align}
&\Pbr{p_i,q_j}=\delta_{i,j},\quad 
\Pbr{q_i,q_j}=\Pbr{p_i,p_j}=0,\\
&\Pbr{p_i,x_0}=\Pbr{q_i,x_0}=\Pbr{p_i,x_1}=\Pbr{q_i,x_1}
=\Pbr{x_0,x_1}=0
\notag
\end{align}
\end{equation}
for $i,j=1,\ldots,n$, we have the same formula as \eqref{Pbr-can} and 
\begin{equation}
\varphi'=\Pbr{h_0,\varphi}+\frac{k}{2}\big(
x_0\frac{\partial \varphi}{\partial x_0}+x_1\frac{\partial \varphi}{\partial x_1}\big) 
\quad\text{for any}\quad \varphi\in\BC(\alpha;f)
\end{equation}
by \eqref{A2n+1-by-Pbr}. 
Let us denote by $H=H(q;p;x)\in\BC(\alpha)[q;p;x^{\pm1}]$ the polynomial in $(q;p)$ 
which represents our Hamiltonian $h_0$ in the coordinates $(q;p;x)$. 
Then we see that the differential system \eqref{A2n+1} is equivalent to the 
Hamiltonian system
\begin{equation}
\frac{dq_i}{dt}=\Pbr{H,q_i}=\frac{\partial H}{\partial p_i},\ 
\frac{dp_i}{dt}=\Pbr{H,p_i}=-\frac{\partial H}{\partial q_i},\ 
\frac{dx_0}{dt}=\frac{k}{2} x_0,\   \frac{dx_1}{dt}=\frac{k}{2} x_1,
\end{equation}
where $i=1,\ldots,n$.  
The Hamiltonian $H$ is given in the form
\begin{equation}\Label{H2n+1}
\begin{align}
H=&
\big(x_0^2-\sum_{i=1}^n q_i \big)\big(\sum_{i=1}^n q_i p_i(\frac{x_1}{x_0}-p_i)\big)
-\sum_{1\le i<j\le n}q_iq_j(p_i-p_j)\big(\frac{x_1}{x_0}+p_i-p_j\big)
\\
&
+2\gamma\sum_{i=1}^n\,q_i p_i
-\frac{x_1}{x_0}\sum_{i=1}^n\beta_i\,q_i
+\sum_{i=1}^n \alpha_{2i}\,p_i
%\notag\\
%&
+(\gamma-\varpi_{2n+1}) x_0x_1 
+\gamma^2,
\notag
\end{align}
\end{equation}
where
\begin{equation}
\beta_i=\sum_{r=1}^i\alpha_{2r-1}\ \ (i=1,\ldots,n),\ \ 
\gamma=\sum_{i=1}^{2n+1}(-1)^{i-1}\varpi_i
=\frac{1}{2}\sum_{r=0}^n \alpha_{2r+1}.
\end{equation}
\medskip
Summarizing the results of this section, we have 

\begin{thm}\Label{thm-Ham} $(1)$ 
The differential system \eqref{A2n} with affine Weyl group symmetry of type $A^{(1)}_{2n}$
is equivalent to the Hamiltonian system 
\begin{equation}
\frac{dq_i}{dt}=\frac{\partial H}{\partial p_i},\quad
\frac{dp_i}{dt}=-\frac{\partial H}{\partial q_i}\quad (i=1,\ldots,n),
\end{equation}
with an auxiliary variable $x$
such that $dx/dt=k$, where the $H=H(q;p;x)$ is the polynomial \eqref{H2n}.
\newline
$(2)$
The differential system \eqref{A2n+1} of type $A^{(1)}_{2n+1}$
is equivalent to the Hamiltonian system 
\begin{equation}
\frac{dq_i}{dt}=\frac{\partial H}{\partial p_i},\ 
\frac{dp_i}{dt}=-\frac{\partial H}{\partial q_i}\quad (i=1,\ldots,n),
\end{equation}
with two auxiliary variables $x_0,x_1$ such that 
$dx_0/dt=k x_0/2$, $dx_1/dt=k x_1/2$,
where $H=H(q;p;x)$ is the polynomial 
\eqref{H2n+1}.
\end{thm}
\noindent
See Appendix for explicit formulas of the Hamiltonians $H$ for small $l$. 

\begin{rem}
When we regard our systems as Hamiltonian systems, it would be more convenient 
to use constant parameters as $x=k t+c$ for $A^{(1)}_{2n}$, and 
$x_0=c_0 e^{kt/2}$, $x_1=c_1 e^{kt/2}$ for $A^{(1)}_{2n+1}$,
rather than the auxiliary variables as in Theorem \ref{thm-Ham}. 
\end{rem}

\section{$\tau$-Functions}

In the following, we introduce a family of Hamiltonians $h_0,\ldots,h_l$ and 
$\tau$-functions $\tau_0,\ldots,\tau_l$ for our system. 
From the Hamiltonian $h_0$ of Section 2,  
we define the other Hamiltonians $h_1,\ldots,h_l$
by the diagram rotation:
\begin{equation}
h_i=\pi(h_{i-1})\quad (i=1,\ldots,l).
\end{equation}
Note that the component of $h_0$ of highest degree in $f_0,\ldots,f_l$ 
is invariant under the diagram rotation.  Hence $h_0,\ldots,h_l$
have a common highest degree component, and have different coefficients 
in lower degrees. 
We also introduce the $\tau$-functions $\tau_0,\ldots,\tau_l$ for our system 
to be the dependent variables such that $h_j=k \big(\log \tau_j \big)'$, namely 
as functions determined from $f_0,\ldots,f_l$ by the linear differential equations
\begin{equation}\Label{tau-h}
k \, \tau_j '= h_j \, \tau_j \quad(j=0,1,\ldots,l).
\end{equation}
In this section, we show that the affine Weyl group symmetry lifts to the level 
of $\tau$-functions as well. 

We first remark that 
our Hamiltonians have some remarkable properties in relation to the action of $W$. 

\begin{prop} \Label{prop-W-on-h}
With respect to the action of the affine Weyl group, the Hamiltonians 
have the following invariance:
\begin{equation}\Label{s-on-h}
\begin{align} 
A^{(1)}_{2n}:&\qquad
s_{i}(h_j)=h_j +\delta_{i,j}\, k\,\frac{\alpha_j}{f_j}\qquad (i,j=0,\ldots,l),
\\
A^{(1)}_{2n+1}:&\qquad
s_{i}(h_j)=h_j +\delta_{i,j}\, k\,\frac{\alpha_j}{f_j}g_j\quad (i,j=0,\ldots,l),
\notag
\end{align}
\end{equation}
where $g_j$ stands for $g_0$ or $g_1$ according as $j\equiv 0$ or $1 \mod 2$. 
\end{prop}
\begin{pf}
We have only to show the case when $j=0$. 
In the case of $A^{(1)}_{2n}$,
we compute $\Delta_i(h_0)$ as follows: 
\allowdisplaybreaks
\begin{eqnarray*}
&& 
\Delta_i(f_{i-1}f_{i}f_{i+1})
+\Delta_i(f_{i-1}f_i)\sum_{K\in {\cal S}_1([i+3,i-2])} f_K
+\Delta_i(f_{i}f_{i+1})\sum_{K\in {\cal S}_1([i+2,i-3])} f_K
\\*
&&
+\Delta_i(f_{i-1})\sum_{K\in {\cal S}_1([i+2,i-2])} f_K
+\Delta_i(f_{i+1})\sum_{K\in {\cal S}_1([i+2,i-2])} f_K
\\*
&&
+\sum_{j=0}^{2n}\Delta_i(\chi(\Gamma\backslash\{j\}))f_j
+s_i(\chi(\Gamma\backslash\{i-1\}))\Delta_i(f_{i-1})
+s_i(\chi(\Gamma\backslash\{i+1\}))\Delta_i(f_{i+1})
\\
&=&
\big(f_{i-1}-f_{i+1}-\frac{\alpha_i}{f_i}\big)
-\sum_{K\in {\cal S}_1([i+3,i-2])} f_K
+\sum_{K\in {\cal S}_1([i+2,i-3])} f_K
\\*
&&
-\frac{1}{f_i}\big(\sum_{K\in {\cal S}_1([i+2,i-2])} f_K\big)
+\frac{1}{f_i}\big(\sum_{K\in {\cal S}_1([i+2,i-2])} f_K\big)
\\*
&&
+\sum_{K\in{\cal S}_1([i+1,i-2])} f_K 
-\sum_{K\in{\cal S}_1([i+2,i-1])} f_K 
\\*
&&
-(\chi(\Gamma\backslash\{i-1\})-\alpha_i)\frac{1}{f_i}
+(\chi(\Gamma\backslash\{i+1\})+\alpha_i)\frac{1}{f_i}
\\*
&=&(\alpha_i+\varpi_{i-1}-2\varpi_i+\varpi_{i+1}) \frac{1}{f_i}
=\delta_{j,0}\,\frac{k}{f_i},
\end{eqnarray*}
which proves \eqref{s-on-h} for $A^{(1)}_{2n}$. 
The case of $A^{(1)}_{2n+1}$ can be verified similarly.  
\end{pf} 

We need also to know how adjacent pairs of the Hamiltonians $h_j$ 
are related.  
We propose some lemmas for this purpose. 

\begin{lem}\Label{lem-chi0} For each $i=0,1,\ldots,l$, one has
\begin{eqnarray}
\pi(\chi([i-1,i]))-\chi([i,i+1]) &=&\begin{cases}
-\frac{1}{l+1}\, k & (i\ne 0),
\\
\frac{l}{l+1}\,k & (i=0).
\end{cases}\\
\pi^{-1}(\chi([i+1,i+2]))-\chi([i,i+1]) &=&\begin{cases}
\frac{1}{l+1}\, k & (i\ne l),
\\
-\frac{l}{l+1}\,k & (i=l).
\end{cases}
\notag
\end{eqnarray}
\end{lem}

\noindent
Let $L$ be a subset of $\Gamma=\{0,1,\ldots,l\}$ with $|L|=2m$ 
and suppose that $L$ is a disjoint union of chains of even nodes. 
We say that $L$ {\em contains} $[i,i+1]$ {\em evenly}, if $[i,i+1]\subset L$ and 
the complement $L\backslash[i,i+1]$ splits into chains of even nodes. 
With this terminology, we have 
\begin{lem}\Label{lem-chi1} 
Under the assumption on $L$ above, one has 
\begin{equation}
\pi(\chi(\pi^{-1} L))-\chi(L)= \frac{l-m+1}{l+1}\,k \quad \text{or}\quad -\frac{m}{l+1}\,k,
\end{equation}
according as $L$ contains $[0,1]$ evenly, or not.  Similarly, 
\begin{equation}
\pi^{-1}(\chi(\pi L))-\chi(L)= -\frac{l-m+1}{l+1}\,k \quad\text{or}\quad \frac{m}{l+1}\,k,
\end{equation}
according as $L$ contains $[l,0]$ evenly, or not. 
\end{lem} 

\noindent
We omit the proof of Lemmas \ref{lem-chi0} and \ref{lem-chi1} since 
they can be proved by direct calculations from the definitions. 

\begin{prop}\Label{prop-h}
$(1)$ In the case of $A^{(1)}_{2n}$ $(l=2n)$, for each $j=0,\ldots,2n$, one has 
\begin{equation}\Label{h-dif-2n}
h_{j+1}-h_j= k\,\sum_{r=1}^{n} f_{j+2r}-\frac{n\,k}{2n+1}\,x,
\end{equation}
where $x=\sum_{i=0}^{2n} f_i$.  Hence
\begin{eqnarray}\Label{h-2n}
-h_{j-1}+2h_j-h_{j+1}&=&k\sum_{r=1}^n \big(f_{j+2r-1}-f_{j+2r}\big),\\
h_{j-1}-h_{j+1}&=&k\, \big(\,f_j-\frac{1}{2n+1} x\,\big). \notag
\end{eqnarray}
$(2)$ In the case of $A^{(1)}_{2n+1}$ $(l=2n+1)$, one has
\begin{eqnarray}\Label{h-dif-2n+1}
h_{j+1}-h_j&=& k\,\sum_{1\le r\le s\le n} f_{j+2r} f_{j+2s+1}\\
&&\quad -\frac{n\,k}{2n+2}\sum_{K\in{\cal S}_2}f_K+
(-1)^j\,\frac{k}{4}\sum_{i=0}^{2n+1}(-1)^{i}\alpha_i,\notag 
\end{eqnarray}
for each $j=0,\ldots,2n+1$. 
Hence
\begin{eqnarray}\Label{h-2n+1}
-h_{j-1}+2h_j-h_{j+1}&=&
k\,\sum_{1\le r\le s\le n}\big(f_{j+2r-1} f_{j+2s}-f_{j+2r} f_{j+2s+1})\\
&& \quad +k\,\big(\,\frac{k}{2}-\sum_{r=0}^n \alpha_{j+2r}\,\big),\notag\\
h_{j-1}-h_{j+1}&=&k \, x_{j+1}\big(\,f_j-\frac{1}{n+1} x_j\,\big), \notag
\end{eqnarray}
where $x_0=\sum_{r=0}^{2n} f_{2r}$ and $x_1=\sum_{r=0}^{2n} f_{2r+1}$,
and $x_i$ stands for $x_0$ or $x_1$ according as $i\equiv 0$ or $1 \mod 2$. 
\end{prop}
\begin{pf}
When $l=2n$, we have 
\begin{eqnarray}
h_0&=&\sum_{K\in{\cal S}_3} f_K+\sum_{i=0}^{2n} \chi([i+1,i-1]) f_i,
\\
h_{1}&=&\sum_{K\in{\cal S}_3} f_K+\sum_{i=0}^{2n} \pi(\chi([i,i-2])) f_i,\notag
\end{eqnarray}
hence
\begin{equation}
h_{1}-h_0=\sum_{i=0}^{2n}\big( \pi(\chi([i,i-2]))-\chi([i+1,i-1])\big)f_i.
\end{equation}
By Lemma \ref{lem-chi1}, we compute
\begin{eqnarray}
h_{1}-h_0&=&-\frac{n\,k}{2n+1}\big(f_0+ \sum_{r=1}^n f_{2r-1})
+\frac{(n+1)\ k}{2n+1}\sum_{r=1}^n f_{2r}
\\ 
&=&k\sum_{r=1}^n f_{2r} - \frac{n\,k}{2n+1}\sum_{i=0}^{2n} f_i, 
\notag
\end{eqnarray}
which gives \eqref{h-dif-2n} for $j=0$.   
Formulas \eqref{h-dif-2n} for $j=1,\ldots,2n$ are obtained by applying 
the diagram rotation $\pi$. 
Formulas \eqref{h-2n} follow directly from \eqref{h-dif-2n}. 
When $l=2n+1$, we have 
\begin{eqnarray}
h_{0}&=&\sum_{K\in{\cal S}_4} f_K
+\sum_{K\in{\cal S}_2}\chi(K^{\text{c}})) f_K+\gamma_1^2,\\
h_{1}&=&\sum_{K\in{\cal S}_4} f_K
+\sum_{K\in{\cal S}_2} \pi(\chi(\pi^{-1}K^{\text{c}})) f_K+\gamma_0^2.\notag
\end{eqnarray}
where $\gamma_1=(\sum_{r=0}^n\alpha_{2r+1})/2$ and 
$\gamma_0=\pi(\gamma_1)=(\sum_{r=0}^n\alpha_{2r})/2$. 
Hence we have 
\begin{equation}
h_{1}-h_{0}=
\sum_{K\in{\cal S}_2}\big(\pi(\chi(\pi^{-1}K^{\text{c}}))- 
\chi(K^{\text{c}})\big)f_K+\frac{k}{2}(\gamma_0-\gamma_1).\notag
\end{equation}
since $\gamma_0+\gamma_1=k/2$. 
The coefficients of $f_K$ are computed by Lemma \ref{lem-chi1} to obtain 
\eqref{h-dif-2n+1} for $j=0$. 
Formulas \eqref{h-dif-2n+1} for the other $j$ are obtained by the diagram rotation, 
and formulas \eqref{h-2n+1} follow directly from \eqref{h-dif-2n+1}. 
\end{pf}

By combining Propositions \ref{prop-W-on-h} and \ref{prop-h}, we obtain 

\begin{prop}\Label{prop-h-motif}
The differential system \eqref{A2n} $($resp. \eqref{A2n+1}$)$ is expressed 
as follows in terms of Hamiltonians $h_0,\ldots,h_l$:
\begin{equation}\Label{h-motif}
k \,\frac{f_j'}{f_j}=s_j(h_j)+ h_j-h_{j-1}-h_{j+1}\quad(j=0,\ldots,l),
\end{equation}
or equivalently, 
\begin{equation}
k \,\frac{f_j'}{f_j}= s_j(h_j)-h_j+\sum_{i=0}^l h_{i}a_{ij}
\end{equation}
where $A=(a_{ij})_{0\le i,j\le l}$ is the generalized Cartan matrix 
of \eqref{GCM}.
\end{prop}
\noindent
Formulas \eqref{h-motif} can be verified by rewriting the right-hand side 
as the sum $(s_j(h_j)-h_j)+(-h_{-j+1}+2 h_j - h_{j+1})$, and then by applying 
Propositions \ref{prop-W-on-h} and \ref{prop-h}. 

\medskip
In what follows, we denoted by $\BC(\alpha;f;\tau)$ the field of rational 
functions in $\tau_0,\ldots,\tau_l$ with coefficients in $\BC(\alpha;f)$. 
We define the structure of differential field of $\BC(\alpha;f;\tau)$
by \eqref{tau-h}. 
We now extend the automorphisms $s_0,\ldots,s_l$ and $\pi$ of $\BC(\alpha;f)$
to $\BC(\alpha;f;\tau)$ by setting 
\begin{equation}
s_i(\tau_j)=\tau_j \ \ (i\ne j),\quad
s_j(\tau_j)=\frac{\tau_{j-1}\ \tau_{j+1}}{\tau_j}\, f_j,
\quad
\pi(\tau_j)=\tau_{j+1} 
\end{equation}
where $i,j=0,\ldots,l$. 
\begin{thm}[\cite{NY3}] The automorphisms $s_0,\ldots, s_l$ and $\pi$ 
of $\BC(\alpha;f;\tau)$ described above define a representation of 
the extended affine Weyl group $\widetilde{W}$. 
\end{thm}

Note that, by Proposition \ref{prop-h}, one obtains an expression 
of $f_j$ by $\tau$ functions.
\begin{prop}\Label{prop-f-in-tau} 
The variables $f_j$  $(j=0,1,\ldots,l)$ are expressed in terms of $\tau$-functions 
as follows:
\begin{eqnarray}
f_j&=&\frac{1}{k}(h_{j-1}-h_{j+1})+\frac{x}{2n+1}\\
&=& \frac{\tau_{j-1}'}{\tau_{j-1}} 
-\frac{\tau_{j+1}'}{\tau_{j+1}}+\frac{x}{2n+1}
\qquad (j=0,\ldots,2n)\notag
\end{eqnarray}
when $l=2n$, and
\begin{eqnarray}
f_j&=&\frac{1}{k \, x_{j+1}}(h_{j-1}-h_{j+1})+
\frac{x_j}{n+1}\\
&=& \frac{1}{x_{j+1}}\big(\frac{\tau_{j-1}'}{\tau_{j-1}} 
-\frac{\tau_{j+1}'}{\tau_{j+1}}\big)+\frac{x_j}{n+1}
\qquad (j=0,\ldots,2n+1)\notag
\end{eqnarray}
when $l=2n+1$. 
\end{prop}
\noindent

Hence we have
 
\begin{prop}
For each $j=0,\ldots,l$, 
the action of $s_j$ on the $\tau$-functions $\tau_j$ is given by the following bilinear operators 
of Hirota type:  
\begin{eqnarray}
s_j(\tau_j)&=&\frac{1}{\tau_j}\big(D_{t}+
\frac{x}{2n+1}\big)\tau_{j-1}\cdot\tau_{j+1}\\
&=&\frac{1}{\tau_j}\big(\tau_{j-1}'\tau_{j+1}-\tau_{j-1}\tau_{j+1}'+
\frac{x}{2n+1}\tau_{j-1}\tau_{j+1}\big)\notag
\end{eqnarray}
when $l=2n$, and  
\begin{eqnarray}
s_j(\tau_j)&=&\frac{1}{\tau_j}\big(
\frac{1}{x_{j+1}}D_t+
\frac{x_j}{n+1}\big) \tau_{j-1}\cdot\tau_{j+1}
\\
&=&\frac{1}{\tau_j}\big(
\frac{1}{x_{j+1}}(\tau_{j-1}'\tau_{j+1}-\tau_{j-1}\tau_{j+1}')+
\frac{x_j}{n+1}\tau_{j-1}\tau_{j+1}
\big)\notag
\end{eqnarray} 
when $l=2n+1$, where $x_j$ stands for $x_0=\sum_{r=0}^n \, f_{2r}$ or 
$x_1=\sum_{r=0}^n \, f_{2r+1}$ according as $j\equiv 0$ or $1 \mod 2$. 
\end{prop}

From Proposition \ref{prop-h-motif}, we obtain 
\begin{thm}\Label{thm-motif}
The action of the extended affine Weyl group $\widetilde{W}$ commutes 
with the derivation of the differential field $\BC(\alpha;f;\tau)$. 
\end{thm}
\begin{pf} Since the equality $s_i(\tau_j)'=s_i(\tau_j')$ is obvious if $i\ne j$, 
we consider the case when $i=j$.  
By Proposition \ref{prop-h-motif}, we compute 
\begin{eqnarray}
\frac{s_j(\tau_j)'}{s_j(\tau_j)}&=&
\frac{\tau_{j-1}'}{\tau_{j-1}}+\frac{\tau_{j+1}'}{\tau_{j+1}}-\frac{\tau_{j}'}{\tau_{j}}
+\frac{f_j'}{f_j}\\
&=& \frac{1}{k}\big(h_{j-1}+h_{j+1}-h_j+k \frac{f_j'}{f_j}\big)\notag\\
&=& \frac{1}{k} s_j(h_j) =s_j(\frac{\tau_j'}{\tau_j})=\frac{s_j(\tau_j')}{s_j(\tau_j)},
\notag
\end{eqnarray}
which implies $s_j(\tau_j)'=s_j(\tau_j')$.  
It is clear that the action of the diagram rotation $\pi$ commutes with $'$. 
\end{pf}

\noindent
Theorem \ref{thm-motif} means that one can lift the B\"acklund transformations 
of our system for the variables $f_j$ 
to the level of $\tau$-functions so that each $\tau_j$ are invariant 
with respect to the subgroup $W_j=\br{s_0,\ldots,s_{j-1},s_{j-1},\ldots,s_l}$ 
of $W$ ($j=0,\ldots,l$),  
and that one has the multiplicative formulas 
\begin{equation}
f_j=\frac{\tau_j\ s_j(\tau_j)}{\tau_{j-1}\tau_{j+1}}\qquad(j=0,\ldots,l)
\end{equation}
for the dependent variables $f_0,\ldots,f_l$ 
in terms of $\tau$-functions. 
We have thus guaranteed that the structure of B\"acklund transformations 
of our differential system is consistent with the general scheme of 
our previous paper \cite{NY3}.  
As a consequence we see that the B\"acklund transformations of our system 
provides the discrete dynamical systems of type $A^{(1)}_l$ in the sense 
of \cite{NY3}.  

\def\thesection{A}
\section{Appendix}

\subsection{Explicit formulas of $f_0'$ for small $l$}

\allowdisplaybreaks
\begin{eqnarray*}
A^{(1)}_2:\quad& f_0'=&f_0(f_1-f_2)+ \alpha_0,\\ 
A^{(1)}_4:\quad& f_0'=&f_0(f_1-f_2+f_3-f_4)+ \alpha_0,\\ 
A^{(1)}_3:\quad& f_0'=&f_0(f_1f_2-f_2f_3)+ 
\big(\frac{k}{2}-\alpha_2\big) f_0+\alpha_0 f_2,\\ 
A^{(1)}_5:\quad& f_0'=&f_0(f_1f_2+f_1f_4+f_3f_4-f_2f_3-f_2f_5-f_4f_5) \\*
&&+\big(\frac{k}{2}-\alpha_2-\alpha_4\big) f_0+\alpha_0 (f_2+f_4).
\end{eqnarray*}

\subsection{Explicit formulas of $h_0$ for small $l$}
\allowdisplaybreaks
\begin{eqnarray*}
A^{(1)}_2:\ & h_0= & 
f_0f_1f_2+\frac{1}{3}(\alpha_1-\alpha_2)f_0+
\frac{1}{3}(\alpha_1+2\alpha_2)f_1-\frac{1}{3}(2\alpha_1+\alpha_2) f_2\\
A^{(1)}_4:\ & h_0= & 
f_0f_1f_2+f_1f_2f_3+f_2f_3f_4+f_3f_4f_0+f_4f_0f_1\\*
&&
+\frac{1}{5}(2\alpha_1-\alpha_2+\alpha_3-2\alpha_4)f_0
+\frac{1}{5}(2\alpha_1+4\alpha_2+\alpha_3+3\alpha_4)f_1\\*
&&
-\frac{1}{5}(3\alpha_1+\alpha_2-\alpha_3+2\alpha_4)f_2
+\frac{1}{5}(2\alpha_1-\alpha_2+\alpha_3+3\alpha_4)f_3\\*
&&
-\frac{1}{5}(3\alpha_1+\alpha_2+4\alpha_3+2\alpha_4)f_4\\
A^{(1)}_3:\ & h_0= & 
f_0f_1f_2f_3
+\frac{1}{4}(\alpha_1+2\alpha_2-\alpha_3)f_0f_1
+\frac{1}{4}(\alpha_1+2\alpha_2+3\alpha_3)f_1f_2\\*
&&
-\frac{1}{4}(3\alpha_1+2\alpha_2+\alpha_3)f_2f_3
+\frac{1}{4}(\alpha_1-2\alpha_2-\alpha_3)f_3f_0
+\frac{1}{4}(\alpha_1+\alpha_3)^2\\
A^{(1)}_5:\ & h_0= & 
f_0f_1f_2f_3+f_1f_2f_3f_4+f_2f_3f_4f_5+f_3f_4f_5f_0+f_4f_5f_0f_1+f_5f_0f_1f_2\\*
&&
+\frac{1}{3}(\alpha_1+2\alpha_2+\alpha_4-\alpha_5)f_0f_1
+\frac{1}{3}(\alpha_1+2\alpha_2+3\alpha_3+\alpha_4+2\alpha_5)f_1f_2\\*
&&
-\frac{1}{3}(2\alpha_1+\alpha_2-\alpha_4+\alpha_5)f_2f_3
+\frac{1}{3}(\alpha_1-\alpha_2+\alpha_4+2\alpha_5)f_3f_4\\*
&&
-\frac{1}{3}(2\alpha_1+\alpha_2+3\alpha_3+2\alpha_4+\alpha_5)f_4f_5
+\frac{1}{3}(\alpha_1-\alpha_2-2\alpha_4-\alpha_5)f_5f_0
\\*
&&
+\frac{1}{3}(\alpha_1-\alpha_2+\alpha_4-\alpha_5)f_0f_3
+\frac{1}{3}(\alpha_1+2\alpha_2+\alpha_4+2\alpha_5)f_1f_4\\*
&&
-\frac{1}{3}(2\alpha_1+\alpha_2+2\alpha_4+\alpha_5)f_2f_5
+\frac{1}{4}(\alpha_1+\alpha_3+\alpha_5)^2
\end{eqnarray*}
 
\subsection{Explicit formulas of $H$ for small $l$}
\allowdisplaybreaks
\begin{eqnarray*}
A^{(1)}_2:\ \ & H=&(x-q_1)q_1p_1-q_1 p_1^2-\alpha_1 q_1+\alpha_2 p_1
+\frac{1}{3}(\alpha_1-\alpha_2) x\\
A^{(1)}_4:\ \ & H=&(x-q_1-q_2)(q_1p_1+q_2p_2)
-q_1 p_1^2-q_2 p_2^2 -q_1(p_1-p_2)q_2 \\*
&&-\alpha_1q_1-(\alpha_1+\alpha_3)q_2+\alpha_2p_1+\alpha_4 p_2+
\frac{1}{5}(2\alpha_1-\alpha_2+\alpha_3-2\alpha_4) x\\
A^{(1)}_3:\ \ 
& H=&  (x_0^2-q_1)p_1q_1(\frac{x_1}{x_0}-p_1)
+(\alpha_1+\alpha_3)q_1p_1 -\alpha_1 \frac{x_1}{x_0}q_1
+\alpha_2 x_0^2p_1
\\*
&&
+\frac{1}{4}(\alpha_1-2\alpha_2-\alpha_3) x_0 x_1
+\frac{1}{4}(\alpha_1+\alpha_3)^2\\
A^{(1)}_5:\ \
& H=&
(x_0^2-q_1-q_2)\big(q_1p_1(\frac{x_1}{x_0}-p_1)+q_2p_2(\frac{x_1}{x_0}-p_2)\big)\\*
&&
-q_1q_2(p_1-p_2)(\frac{x_1}{x_0}+p_1-p_2)
%\\*
%&&
+(\alpha_1+\alpha_3+\alpha_5)(q_1p_1+q_2 p_2)\\*
&&
-\frac{x_1}{x_0}\big(\alpha_1  q_1+(\alpha_1+\alpha_3)q_2\big)
+x_0^2(\alpha_2 p_1+\alpha_4 p_2)\\*
&&
+\frac{1}{3}(\alpha_1-\alpha_2-2\alpha_4-\alpha_5) x_0 x_1 
+\frac{1}{4}(\alpha_1+\alpha_3+\alpha_5)^2
\end{eqnarray*}
%%%%%%%%%%%%%%%%%%%%%%%%%%%%%%%%%%%%%%%%%%%%%%%%%%%

%%%%%%%%%%%%%%%%%%%%%%%%%%%%%%%%%%%%%%%%%%%%%%%%%%%

\end{document}